\documentclass{amsart}

\usepackage{latexsym}
\usepackage{amssymb}

\newtheorem{THM}{Theorem}

\newtheorem{Cor}[THM]{Corrolary}

\renewcommand{\Re}{\mathbb{R}}

\renewcommand{\phi}{\varphi}
\renewcommand{\S}{\mathbb{S}}

\newcommand{\ra}{\rightarrow}

\newcommand{\lp}{\left(}
\newcommand{\rp}{\right)}
\newcommand{\lb}{\left[}
\newcommand{\rb}{\right]}

\newcommand{\floor}[1]{\lfloor#1\rfloor}

\DeclareMathOperator{\dist}{dist}
\DeclareMathOperator{\sign}{sgn}

\begin{document}

\title{Brownian motion and the parabolicity of minimal graphs}
\author{Robert W.\ Neel}
\address{Department of Mathematics, Columbia University, New York, NY}
\begin{abstract}We prove that minimal graphs (other than planes) are parabolic in the sense that any bounded harmonic function is determined by its boundary values.  The proof relies on using the coupling introduced in~\cite{MyPaper} to show that Brownian motion on such a minimal graph almost surely strikes the boundary in finite time.
\end{abstract}
\thanks{The author gratefully acknowledges support from an NSF Postdoctoral Research Fellowship.}
\email{neel@math.columbia.edu}
\date{October 3, 2008}

\maketitle

\section{Introduction}

By a \emph{minimal graph}, we mean a complete minimal surface-with-boundary, the interior of which can be written as a graph over some open subset of the plane.  (We assume that our minimal graphs are connected.)  A famous result of Bernstein (see~\cite{Osserman}) states that the only minimal graphs over the entire plane are planes (that is, the graphs of affine functions).  In this paper, we prove that all other minimal graphs (aside from the planes just mentioned) are parabolic, in the sense that any bounded harmonic function is determined by its boundary values.  In particular, this means that the boundary of such a surface is non-empty.  We obtain this result as a corollary of proving that Brownian motion on any non-planar minimal graph almost surely has finite lifetime, which means that it almost surely strikes the boundary in finite time.  Our proof of this is essentially an application the coupling introduced in~\cite{MyPaper}.

With this result, we settle a conjecture (in the affirmative) due to Meeks that any minimal graph over a proper subdomain of the plane is parabolic.  Related questions have been studied by L\'{o}pez and P\'{e}rez~\cite{LopezPerez}, who prove that a non-flat properly immersed minimal surface-with-boundary that lies above a sublinear graph and that has Gauss map contained in an open hyperbolic subset of the sphere is parabolic.

\section{A description of the coupling}

As mentioned, our proof leans heavily on the results of~\cite{MyPaper}.  We summarize the relevant results, specialized to the case of a single minimal surface-with-boundary, here; the interested reader is encouraged to consult the original paper for more details.

Let $M$ be a minimal surface-with-boundary.  Then, for any point $(x_0,y_0)$ in the interior of $M\times M$, a coupled Brownian motion is a continuous stochastic process $(x_t,y_t)$, starting at $(x_0,y_0)$, such that the $x_t$ and $y_t$ marginals are both Brownian motions on $M$, defined until the first time $x_t=y_t$ or either of the marginals strikes the boundary of $M$.  In order to describe the particular coupling we will use, we need to introduce some processes associated to any coupled Brownian motion.  First, let $r_t=\dist_{\Re^3}(x_t,y_t)$ be the (extrinsic) distance between the particles.  Next, let $m:M\ra \S^2$ be the Gauss map of $M$.  The coupling will be governed by the relative positions of three unit vectors (which we think of as points in $\S^2$), namely $m(x_t)$, $m(y_t)$, and $\alpha(x_t,y_t)=(x_t-y_t)/|x_t-y_t|$.  We will refer to the position of these three unit vectors as the configuration of the system.

In what follows, we will use $W_{z}$ to denote a Brownian motion (on the real line) with time parameter $z$.  The particular Brownian motion may vary from use to use, and we will use $\widetilde{W_{z}}$ to denote a second such Brownian motion (which need not be independent of the first).  In~\cite{MyPaper}, it was proved that there exists a coupled Brownian motion on $M$, starting from any $(x_0,y_0)$, such that the semi-martingale decomposition of $r_t$ can be written as
\[
dr_t = \sqrt{f} \,dW_t + \frac{g}{2r_t} \,dt
\]
where $f$ and $g$ are non-negative functions of $m(x_t)$, $m(y_t)$, and $\alpha(x_t,y_t)$ (with one caveat, described below) satisfying several additional properties which we now describe.  (That $f$ and $g$ depend only on the configuration was not required of the coupled Brownian motions considered in~\cite{MyPaper}, but it's easy to see that this can be arranged.)  First, the inequality $f\geq g$ holds everywhere, and we take a moment to explain the significance of this.  If we introduce the time-change $\tau(t) = \int_0^t f ds$, then the martingale part of $r_{\tau}$ is a Brownian motion (we somewhat abuse notation and use $r_{\tau}$ to denote what should be $r_{t(\tau)})$ and the inequality $f\geq g$ means that $r_{\tau}$ is dominated by a two-dimensional Bessel process.  Further, if the inequality is strict at any instant, then the domination is also instantaneously strict.

Second, we can characterize the configurations where $f=g$.  They all occur when $m(x_t)$, $m(y_t)$, and $\alpha(x_y,y_t)$ lie on the same great circle.  Let $\tilde{S}$ be the subset of $(\S^2)^3$ such that all three points lie on a great circle in $\S^2$.  On any simply-connected subset of $\tilde{S}$, we can define a coordinate for the great circle by the signed distance from $\alpha(x_t,y_t)$ (obviously, there are two possible choices of coordinate related by a change of sign).  Then let $\theta$ be the coordinate at $m(x_t)$ and let $\phi$ be the coordinate at $m(y_t)$.  We would like $\theta$ and $\phi$ to be continuous, and thus we will think of them as being defined only up to multiples of $2\pi$.  Since we will only be concerned with the values of trig functions at $\theta$ and $\phi$, defining them only up to $2\pi$ won't cause any problems.  We have that, on this component of $\tilde{S}$,
\[\begin{split}
f & = \lp\sin\theta -A\sin\phi\rp^2 \quad\text{and}\quad 
g = \lp\cos\theta -A\cos\phi\rp^2 \\
\text{where} & \quad A=\sign\lp\cos(\theta+\phi)\rp .
\end{split}\]
Here the choice of whether $A$ is $-1$ or $1$ at a point $(x,y)$ with $\cos(\theta+\phi)=0$ depends on the geometry of $M\times M$ near $(x,y)$, and this is the one caveat to our statement that the evolution of $r_t$ (instantaneously) depends only on $r_t$ and the relative positions of $m(x_t)$, $m(y_t)$, and $\alpha(x_t,y_t)$.  Nonetheless, the determination of $A$ at such points isn't relevant to this paper.  Continuing, we see that $f=g=0$ when $\theta=\phi$ and $A=1$.  Also, when $\cos(\theta+\phi)=0$, we see that there are two possibilities for $f$ and $g$, corresponding to $A=1$ and $A=-1$, and for each possibility we have $f=g>0$, with one exception.  Namely, when $\cos(\theta+\phi)=0$ and $\theta=\phi$, we have that $A=1$ gives $f=g=0$ and $A=-1$ gives $f=g=2$.

The above gives us good control of the evolution of $r_t$ for configurations in $\tilde{S}$.  Our next goal is to extend this to a larger set of configurations.  We introduce a sequence of subsets of $(\S^2)^3$ which we will use throughout the remainder of the paper.  Let $S_1$ be the set where $m(x_t)$ and $m(y_t)$ are at least some fixed (small) positive distance $c_1$ apart in the $\S^2$ metric.  (We won't specify a value for $c_1$; rather, we will assume that it satisfies various properties as we go.)  Let $S_4$ be the subset of $\tilde{S}$ where $\cos(\theta+\phi)=0$ and where $m(x_t)$ and $m(y_t)$ are at least distance $2c_1$ apart.  Let $S_3$ be a closed neighborhood of $S_4$.  Let $S_2$ be an closed set containing $S_3$ on which $|\cos(\theta+\phi)|\leq c_2$ for some positive constant $c_2$, and such that the boundary of $S_2$ is positive distance from the boundary of $S_3$, in the product metric on $(\S^2)^3$.  Further, we assume that $S_2$ and $S_3$ are chosen so that
\[
S_4 \subset S_3 \subset S_2 \subset S_1
\]
and so that the boundaries of all four sets are a positive distance from each other.  Much as we do for $c_1$, we will feel free to adjust $S_3$ and $S_2$ as we go.

We extend $\theta$ and $\phi$ from $S_4$ to be smooth functions on all of $S_2$.  (Note that this extension is not the same as used in~\cite{MyPaper}, but this is a more convenient definition of $\theta$ and $\phi$ for our present purpose.)  On $S_3$ we have that $f-g\geq c_3|\cos(\theta+\phi)|$ for some positive constant $c_3$, while on $S_1\setminus S_3$, we have that $f-g\geq c_4$ for some positive constant $c_4$.

Now that we have good control of $f-g$ in terms of $\theta+\phi$, we need to understand how the configuration evolves.  Again, we will simply cite facts from~\cite{MyPaper}.  Both $m(x_t)$ and $m(y_t)$ are time-changed spherical Brownian motions, with the time-changes given by the integrals of $-K$ along the respective paths on $M$.  Further, $\alpha(x_t,y_t)$ is a semi-martingale on $\S^2$, such that both the quadratic variation and the drift grow as a rate bounded from above by a multiple of $1/r_t^2$.  This gives us control over how quickly the configuration can move from one region of $(\S^2)^3$ to another in terms of $K$ and $r_t$.  In the other direction, we will need to know that the configuration doesn't spend too much time in $S_3$.  Recall the time-change $\tau(t)=\int_0^t f \,dt$.  Since $f$ is bounded above and below by positive constants on $S_1$, so is $d\tau/dt$.  In particular, if we assume that the configuration of our coupled Brownian motion never leaves $S_1$, then a process has infinite lifetime in the original time $t$ if and only if it also has infinite lifetime in the $\tau$ time-scale.  Also, all of our earlier estimates on the evolution of processes in the $t$ time-scale hold, up to constants, in the $\tau$ time-scale, assuming the configuration stays in $S_1$.  We have the following semi-martingale decomposition for $(\theta+\phi)_{\tau}$, valid on $S_2$,
\[
d(\theta+\phi)_{\tau} = \frac{2+\epsilon_1}{r_{\tau}}\, dW_{\tau} + a(x_{\tau},y_{\tau})\, d\widetilde{W}_{\tau} +\lp A\frac{2+\epsilon_2}{r^2_{\tau}} +b(x_{\tau},y_{\tau})\rp\, d\tau 
\]
where $A$ takes only the values $-1$ and 1 and extends our earlier definition of $A$ on $\tilde{S}$ to $S_2$.  In particular, $A$ is $-1$ on some open subset of $S_2$, 1 on another open subset, and its value on their common boundary depends on the geometry of $M\times M$.  Here $\epsilon_1$ and $\epsilon_2$ are bounded functions on $S_2$, with common bound $\epsilon>0$ that can be made arbitrarily small by shrinking $S_2$.  Further, $a$ and $b$ are functions depending on the Gauss map at $x$ and $y$ with the property that both $a^2$ and $b$ are bounded pointwise by a multiple of the sum of the absolute values of the Gauss curvatures at $x$ and $y$, with the consequence that $\int a^2 \,d\tau$ and $\int b\, d\tau$ along any path (in $M\times M$) are bounded by a multiple of the sum of the integrals of the absolute values of the Gauss curvatures along the $x_t$ and $y_t$ marginals.

\section{The coupling on minimal graphs}

We now have enough information about the behavior of the coupled Brownian motion to prove our main theorem.  Note that by \emph{embedded} we mean injectively immersed.

\begin{THM}
Let $M$ be an embedded, complete minimal surface-with-boundary which is not a plane.  Further assume that, outside of a compact subset of $M$, the Gauss map of $M$ is contained in an open, hyperbolic subset of the sphere.  Then Brownian motion on $M$ almost surely hits the boundary in finite time, and thus $M$ is parabolic. 
\end{THM}

\emph{Proof:}  The theorem is clear when $M$ is (isometric to) a subset of the plane.  Thus we assume that $M$ is not flat, and we have that the Gauss curvature, $K\leq 0$, has only isolated zeroes.

Let $B_t$ be a Brownian motion on $M$, started at a point $p$ in the interior.  Let $\sigma$ be the first hitting time of the boundary; we will always stop the Brownian motion at the boundary.  We wish to prove that $\sigma$ is almost surely finite.  The induced process on the Gauss sphere, $m(B_{t\wedge\sigma})$, is a time-changed Brownian motion on $\S^2$, with the time-change given by $u(t\wedge\sigma)=\int_0^{t\wedge\sigma} -K\circ B_v \,dv$, the integral of the absolute value of the Gauss curvature along the path.  Our assumptions on the image of $m$ imply that $u(t\wedge\sigma)$ almost surely converges to a finite limit as $t\ra\infty$ and thus that the process $m(B_{t\wedge\sigma})$ almost surely converges.  Further, using that $K$ has only isolated zeroes, we see that the limiting distribution of the process on the sphere, which we denote $m(B_{\infty\wedge\sigma})$ (despite the fact that $B_{\infty\wedge\sigma}$ is not itself well-defined a priori), does not charge any points.

We now suppose that $\sigma$ is not almost surely finite and derive a contradiction.  Under this assumption, we can find starting points for Brownian motion on $M$ for which the probability of $\{\sigma=\infty\}$ is arbitrarily close to 1.  Further, because $m(B_{\infty\wedge\sigma})$ does not charge points, among such starting points we can find two points $x_0$ and $y_0$ with the following further property: there exist subsets of $\S^2$, $\Gamma_1$ and $\Gamma_2$, such that the four sets $\pm\Gamma_1$ and $\pm\Gamma_2$ are positive distance apart and such that Brownian motion started at $x_0$ has normal vector that stays in $\Gamma_1$ for its entire lifetime with probability arbitrarily close to 1, and similarly for $y_0$ and $\Gamma_2$.  It follows that we can find points $x_0$ and $y_0$, with corresponding sets $\Gamma_1$ and $\Gamma_2$, such that the coupled Brownian motion described above, started from $(x_0,y_0)$, has probability at least $1/2$ of never hitting the boundary and of $(m(x_{t\wedge\sigma}), m(y_{t\wedge\sigma}))$ staying in $\Gamma_1\times\Gamma_2$ for all time.  Further, because $M$ is embedded, note that $x_t$ and $y_t$ cannot couple (that is, $r_t$ cannot strike zero) when $(m(x_{t\wedge\sigma}), m(y_{t\wedge\sigma}))$ is in $\Gamma_1\times\Gamma_2$.

Next, we wish to study our adequate coupling, as described in the previous section, started at $(x_0,y_0)$.  Let $\tilde{\sigma}$ be the first time the process either hits the boundary or has $(m(x_t),m(y_t))$ exit $\Gamma_1\times\Gamma_2$.  Let $\Omega$ be the set of paths (under the adequate coupling started at $(x_0, y_0)$) for which $\tilde{\sigma}=\infty$.  Then paths in $\Omega$ have infinite lifetime, since they also never couple, and the probability of $\Omega$ is at least $1/2$.  We can choose a set $S_1$, as in the previous section, such that all paths in $\Omega$ have configurations which stay in $S_1$ for all time.  Since $f$ is bounded from below on $S_1$ and the paths in $\Omega$ have infinite lifetime, it follows that $r_{\tau}$ almost surely has infinite quadratic variation on $\Omega$.  Because $r_{\tau}$ is dominated (from above) by a two-dimensional Bessel process and is bounded from below by a Brownian motion, we know that infinite quadratic variation means that $r_{\tau}$ almost surely hits every positive level infinitely often.

The contradiction we're working toward is to show that $r_{\tau}$ almost surely does not hit every positive level infinitely often on $\Omega$.  To see this, we introduce the change of coordinate $\rho=\log r$ and the time-change $s=\int (1/r^2_{\tau})\, d\tau$.  Then we have
\begin{equation}\label{Eqn:RhoEvo}
d\rho_s = dW_s - \frac{1-g/f}{2}\, ds .
\end{equation}
Further, we have subsets $S_2$, $S_3$, and $S_4$ of $S_1$ as described in the previous section such that
\[
d(\theta+\phi)_s = (2+\epsilon_1)\, dW_s + a r\, d\widetilde{W}_s + \lb A(2+\epsilon_2) +br^2\rb \, ds
\]
on $S_2$.  Here we note that a simple computation, using that the integrals of the Gauss curvature along the two marginals are almost surely finite, shows that the integrals $\int a^2r^2\, ds$ and $\int br^2 \,ds$ are almost surely finite.

Note that, because $f$ is bounded from below on $S_1$, our inequalities for $g/f$ imply similar inequalities for $1-g/f$.  In particular, we see from Equation~\eqref{Eqn:RhoEvo} that $\rho_s$ is a Brownian motion with non-positive drift and that this drift is bounded from above by a negative constant on $S_1\setminus S_3$.  We wish to show that this drift, in some average sense, spends significant time away from zero.  We partition $s$-time into the intervals $I_n= [n-1, n)$ for all non-negative integers $n$ and consider the behavior of $(\theta+\phi)_s$ over these intervals, since this will allow us to estimate the drift of $\rho_s$.

We proceed by a series of special cases.  First, assume that $a$ and $b$ are identically zero, and consider the process until $\tilde{\sigma}$.  We have $d(\theta+\phi)_s = (2+\epsilon_1)\, dW_s + A(2+\epsilon_2)\, ds$.  We wish to prove that the drift (of $\rho_s$) has a fixed probability of being less than a negative constant over $I_n$.  This follows from showing that, during $I_n$, the configuration has a positive probability of spending some positive amount of time outside of $S_3$.  By a change of measure using Girsanov's theorem, we can eliminate the drift of $(\theta+\phi)_s$ so that it becomes a time-changed Brownian motion, with time change (relative to $s$) bounded below by $(2-\epsilon)^2$ and above by $(2+\epsilon)^2$.  Recall that $|\cos(\theta+\phi)|<c_2$ on $S_2$.  Thus we see that the configuration has a positive probability of leaving $S_2$ within time $\delta_1>0$, where $\delta_1$ can be made as small as we wish by choosing $c_2$ to be small.  Next, note that if the configuration is in $S_1\setminus S_2$, it has a positive probability of taking at least time $\delta_2>0$ to return to $S_3$.  This follows from the fact that the boundaries of $S_2$ and $S_3$ are positive distance apart along with our earlier observations about the evolution of the configuration, adjusted for the change to $s$-time.  Finally, we need some understanding of how these estimates relate to one another.  Our coupled Brownian motion and its associated configuration process are not necessarily Markov.  Nonetheless, these estimates (on how quickly the configuration exits $S_2$ and how quickly it returns to $S_3$) hold for any process, as long as the coefficients of its semi-martingale decomposition obey the right estimates, and are thus independent of the past.  (Intuitively, these estimates are the result of implicitly comparing our processes with the ``worst-case'' processes, and these comparison processes are Markov.)  We conclude that, after choosing our various constants and sets appropriately (by which we mean the $c_i$, the $\delta_i$, and the $S_i$), there is a positive probability, call it $\gamma$, that the drift of $\rho_s$ will be less that $-\delta$, for some positive $\delta$, on any interval $I_n$.  Further, by the same reasoning as above, the estimate holds on each interval independently of the others.  Thus, there is a countable sequence $X_n$ of independent Bernoulli random variables, each of which is equal to 1 with probability $\gamma$ and 0 with probability $1-\gamma$, such that
\begin{equation}\label{Eqn:Drift}
\int_{0}^{s_0\wedge\tilde{\sigma}} -\frac{1-g/f}{2}\, ds \leq -\delta\sum_{n=1}^{\floor{s_0\wedge\tilde{\sigma}}} X_n
\end{equation}
where the integral is along paths and $\floor{z}$ is the largest integer less than or equal to $z$.

We now wish to lift the requirement that $a$ and $b$ are identically zero.  Instead, we assume that $\int_{I_n} a^2r^2\, ds$ and $\int_{I_n} br^2\, ds$ are less than some small, positive $\tilde{\epsilon}$ over any interval $I_n$ with $n\leq\tilde{\sigma}$.  If $\tilde{\epsilon}$ is small enough, then with probability arbitrarily close to one, the contribution of the terms involving $a$ and $b$ to the evolution of the configuration can be made arbitrarily small.  Thus by assuming $\tilde{\epsilon}$ is sufficiently small, and perhaps making $\gamma$ and $\delta$ slightly smaller and further adjusting our constants $\delta_i$ and $c_i$ and our sets $S_i$, the estimate in Equation~\eqref{Eqn:Drift} still holds.  Applying the strong law of large numbers to the sequence $X_1, X_2,\ldots$, we see that, almost surely for paths in $\Omega$, $\int -(1-g/f)/2 \, ds$ goes to negative infinity linearly.  That is, there exists some positive constant $c$ such that
\begin{equation}\label{Eqn:DriftTwo}
\int_0^{s_0} -\frac{1-g/f}{2}\, ds +cs_0 \ra -\infty \quad\text{as $s_0\ra \infty$}
\end{equation}
for almost every path in $\Omega$.

Finally, we come to the general case, where we make no additional assumptions on $a$ and $b$.  As mentioned above, we know that the integrals $\int a^2r^2\, ds$ and $\int br^2 \,ds$ are almost surely finite.  Thus, for any choice of $\tilde{\epsilon}$, we have that almost every path in $\Omega$ has $\int_{I_n} a^2r^2\, ds$ and $\int_{I_n} br^2 \,ds$ bounded by $\tilde{\epsilon}$ for all but finitely many $n$.  Throwing away finitely many intervals doesn't effect the asymptotic behavior of the drift, so we see that Equation~\eqref{Eqn:DriftTwo} holds in general.

The law of the iterated logarithm implies that $|W_{s_0}|$ almost surely grows sublinearly, in the sense that $W_{s_0}/s_0\ra 0$ almost surely as $s_0\ra\infty$.  It follows from Equation~\eqref{Eqn:RhoEvo} and Equation~\eqref{Eqn:DriftTwo} that $\rho_s\ra-\infty$ as $s\ra\infty$ almost surely on $\Omega$.  This means that for these paths $\rho_s$ has a last time above $0$, which in turn means that $r_{\tau}$ has a last time above $1$.  This contradicts our earlier observation that almost every path in $\Omega$ has infinitely many excursions above any level, and we conclude that Brownian motion on $M$ almost surely strikes the boundary in finite time.  

The last claim of the theorem, that this implies parabolicity, follows from the standard representation of bounded harmonic functions in terms of Brownian motion evaluated at stopping times.  $\Box$

Since any graph is embedded and has its Gauss map restricted to a hemisphere, any minimal graph, other than a plane, satisfies the assumptions of the above theorem.  Thus we've succeeded in proving the following result.

\begin{Cor}
Any minimal graph which is not a plane is parabolic.
\end{Cor}

\section{Acknowledgements}

The author would like to thank Bill Meeks for introducing him to this problem and suggesting that Brownian motion might be a useful technique.

\bibliographystyle{amsplain}

\providecommand{\bysame}{\leavevmode\hbox to3em{\hrulefill}\thinspace}
\providecommand{\MR}{\relax\ifhmode\unskip\space\fi MR }
\providecommand{\MRhref}[2]{%
  \href{http://www.ams.org/mathscinet-getitem?mr=#1}{#2}
}
\providecommand{\href}[2]{#2}

\end{document}